\documentclass[letterpaper, 10 pt, conference]{ieeeconf}  \IEEEoverridecommandlockouts                              
\usepackage{amsmath}
\usepackage{amssymb}
\usepackage{amsfonts}
\usepackage{graphicx}
\usepackage[hidelinks]{hyperref}
\usepackage{color}
\usepackage{algorithm}
\usepackage{cite}
\usepackage{algpseudocode}
\usepackage{booktabs}
\usepackage{multirow}
\usepackage{verbatim}

\allowdisplaybreaks[4]
\usepackage{theorem}
\newtheorem{remark}{Remark}

\title{\LARGE \bf
PIQP: A Proximal Interior-Point Quadratic Programming Solver
}

\author{Roland Schwan, Yuning Jiang, Daniel Kuhn, and Colin N. Jones
\thanks{Open source: \url{https://github.com/PREDICT-EPFL/piqp}}%
\thanks{This work was supported by the Swiss National Science Foundation under the NCCR Automation (grant agreement 51NF40\_180545).
Roland Schwan, Yuning Jiang, and Colin N. Jones are with the Automatic Control Lab, EPFL, Switzerland. Roland Schwan and Daniel Kuhn are with the Risk Analytics and Optimization Chair, EPFL, Switzerland. \tt{\{roland.schwan, yuning.jiang, daniel.kuhn, colin.jones\}@epfl.ch}}
}%

\begin{document}

\maketitle
\setlength\abovedisplayskip{4pt}
\setlength\belowdisplayskip{4pt}

\begin{abstract}
This paper presents PIQP, a high-performance toolkit for solving generic sparse quadratic programs (QP). Combining an infeasible Interior Point Method (IPM) with the Proximal Method of Multipliers (PMM), the algorithm can handle ill-conditioned convex QP problems without the need for linear independence of the constraints. The open-source implementation is written in C++ with interfaces to C, Python, Matlab, and R leveraging the Eigen3 library. The method uses a pivoting-free factorization routine and allocation-free updates of the problem data, making the solver suitable for embedded applications. The solver is evaluated on the Maros-Mészáros problem set and optimal control problems, demonstrating state-of-the-art performance for both small and large-scale problems, outperforming commercial and open-source solvers.
\end{abstract}

\section{Introduction}
Convex quadratic programs are fundamental in many areas of applied mathematics and engineering. They are utilized in various applications, including portfolio optimization~\cite{rubinstein2002markowitz}, optimal control~\cite{rawlings2017model}, state estimation~\cite{allgower1999nonlinear}, and geometry processing~\cite{zhu2018blended}. Furthermore, QPs are a crucial building block of powerful optimization techniques, such as sequential quadratic programming~\cite{boggs1995sequential} for nonlinear programming and branch-and-bound methods~\cite{fletcher1998numerical} for mixed integer quadratic programming. Due to their widespread use, the demand for efficient QP solvers that are both fast and reliable has increased, driven by emerging applications in areas such as optimal control, embedded systems, and signal processing.

In recent decades, significant research efforts have focused on developing efficient QP solvers. Numerous algorithms have been proposed, such as the classical active-set~\cite{bemporad2015quadratic,cimini2019complexity,ferreau2014qpoases} and interior-point methods~\cite{wright1997,domahidi2013,frison2020hpipm}, first-order~\cite{goldstein2014fast,patrinos2013accelerated,patrinos2015dual,stellato2020} and second-order Newton-type methods~\cite{frasch2015parallel,patrinos2011global}, to design QP solvers that achieve high computational efficiency and scalability. As one of the most classical QP approaches, active-set based solvers, such as \texttt{qpOASES}~\cite{ferreau2014qpoases}, are known for their speed in solving small to medium-sized problems and the ability to warm-start using an estimate of the active constraint set. However, they have limited scalability, struggle to exploit sparsity, and are not robust to early termination. Compared to the active-set method, interior-point based solvers, such as \texttt{ECOS} and \texttt{QPSWIFT}, can be fast for large-scale problems and robust to early termination. They can solve sparse large-scale problems efficiently using advanced linear algebra routines but are difficult to warm-start.

In contrast to classical approaches, specifically, interior-point methods, first-order solvers - including operator-splitting-based solvers like~\texttt{OSQP}~\cite{stellato2020} and \texttt{PROXQP}~\cite{bambade2022} - can easily be warm-started, offer simplicity and tight complexity bounds. However, their convergence rates are slower than other methods. In contrast to this, the Newton-type solvers, such as the dual-Newton solver \texttt{qpDUNES}~\cite{frasch2015parallel}, can exploit sparsity and warm-start, but require restrictive assumptions, such as the QP being strongly convex and satisfying the linear independence constraint qualification (LICQ). These assumptions reduce their applicability and may cause robustness issues. Although the primal-dual Newton-type methods, such as \texttt{FBRS}~\cite{liao2018regularized}, can relax the strong convexity requirement to the second-order sufficient condition (SOSC). However, they still require the LICQ.

Recently, hybrid methods that combine first-order and active set/interior-point methods have been the focus of increasing research attention. An example of such a method is \texttt{QPNNLS}~\cite{bemporad2017numerically}. Using the proximal point algorithm, a sequence of regularized QP problems is solved which converges to the solution of the original problem. The regularized QP problems are strictly convex and are solved using a non-negative least-square based active set method. Since solving QP subproblems is computationally expensive, \texttt{QPNNLS} relies heavily on warm-starting to reduce computational cost. Another related method is \texttt{QPALM}~\cite{hermans2022qpalm}, which is based on the Augmented Lagrangian Method. Additionally, \texttt{FBstab}~\cite{liao2020fbstab}, a proximally stabilized Fischer-Burmeister method-based solver, employs a primal-dual version of the proximal point algorithm, in which the proximal subproblems are solved by using a Newton-type method.

In this paper, we present a software contribution, a hybrid approach based QP solver, called \texttt{PIQP}. The underlying algorithm implemented in \texttt{PIQP} follows the framework proposed in~\cite{pougkakiotis2021}, which combines the interior-point method and the proximal method of multipliers (PMM). In particular, it uses one-iteration of Mohetra's predictor-corrector method to deal with the proximal subproblem combining the dual gradient update.

Section~\ref{sec:pre} introduces some preliminaries, including the problem formulation and the main idea of PMM. The practical algorithm implemented in~\texttt{PIQP} is presented in Section~\ref{sec:pipm}. Section~\ref{sec:imp} elaborates the numerical implementation details of \texttt{PIQP}. We demonstrate the effectiveness of our solver in Section~\ref{sec:num}, in which it is compared against five existing state-of-art approaches, including two commercial solvers~\texttt{Gurobi} and~\texttt{Mosek}, three open-source solvers~\texttt{OSQP},~\texttt{SCS}, and \texttt{PROXQP}. All solvers are evaluated on the Maros-Mészáros benchmark problems~\cite{maros1999}.%

\textit{Notations:} we denote the set of real numbers by $\mathbb{R}$, the set of $n$-dimensional real-valued vectors by $\mathbb{R}^n$, and the set of $n \times m$-dimensional real-valued matrices by $\mathbb{R}^{n \times m}$. Moreover, we denote the subspace of symmetric matrices in $\mathbb{R}^{n \times n}$ by $\mathbb{S}^n$ and the cone of positive semi-definite matrices by $\mathbb{S}_+^n$. We denote the identity matrix as $I_n \in \mathbb{R}^{n \times n}$ and a vector filled with ones as $\mathbf{1}_n \in \mathbb{R}^n$. The $\circ$ operator indicates the element-wise multiplication of two vectors.

\section{Preliminaries}
\label{sec:pre}

This section defines the problem formulation considered in~\texttt{PIQP} and briefly reviews the proximal method of multipliers, which will be used later to design the underlying algorithm in \texttt{PIQP}.

\subsection{Problem Formulation}

\texttt{PIQP} considers quadratic programs in the form
\begin{subequations} 
\label{eq:qp_primal}
\begin{align}\label{eq:primalobj}
\min_{x} \quad & \frac{1}{2} x^\top P x + c^\top x \\\label{eq:primaleq}
\text {s.t.}\quad & Ax=b, \\\label{eq:primalineq}
& Gx \leq h,
\end{align}
\end{subequations}
with primal decision variables $x \in \mathbb{R}^n$, matrices $P\in \mathbb{S}_+^n$, $A \in \mathbb{R}^{p \times n}$,  $G \in \mathbb{R}^{m \times n}$, and vectors $c \in \mathbb{R}^n$, $b \in \mathbb{R}^p$, and $h \in \mathbb{R}^m$.
To design a practical numerical solver, it is convenient to rewrite~\eqref{eq:qp_primal} as the equivalent standard form problem
\begin{subequations} 
\label{eq:qp_primal_slack}
\begin{align}
\min_{x,s} \quad & \frac{1}{2} x^\top P x + c^\top x \\\label{eq:liftedeq}
\text{s.t.}\quad & Ax=b, \\\label{eq:liftedineq}
& Gx - h + s = 0, \\\label{eq:liftedslack}
& s \geq 0,
\end{align}
\end{subequations}
where slack variables $s \in \mathbb{R}^m$ are introduced to lift the affine inequality~\eqref{eq:primalineq} into the equality~\eqref{eq:liftedineq}. Although introducing slack variables $s$ may compromise strong convexity when $P\succ 0$, it has distinct computational benefits. For example, it makes it particularly easy to project onto the feasible set in the context of operator splitting-based approaches, this reformulation results in an easily computable projection operator; another example is the well-known interior-point method to be discussed in the next section. In the following, we mainly work with formulation \eqref{eq:qp_primal_slack}.

\subsection{Proximal method of multipliers}
The augmented Lagrangian of problem \eqref{eq:qp_primal} is defined as
\begin{align}\notag
\label{eq:AL}
\mathcal{L}^\text{ALM}_{\delta}(x,s;& \lambda,\nu) := \frac{1}{2} x^\top P x + c^\top x\\\notag
&+\lambda^\top (A x-b)+\frac{1}{2 \delta}\|A x-b\|_2^2 \\
&+\nu^\top (Gx-h+s)+\frac{1}{2 \delta}\|Gx-h+s\|_2^2,
\end{align}
where variables $\lambda$ and $\nu$ are the Lagrangian multipliers of equality constraints~\eqref{eq:liftedeq} and~\eqref{eq:liftedineq}, respectively, and $\delta>0$ is a penalty parameter. We then introduce the iterations of the proximal method of multipliers originally proposed in~\cite{rockafellar1976} as follows:
\begin{subequations}
\label{eq:proximalM}
\begin{align}\label{eq:PMM_sub_prob}
(x^+,s^+) & \in \underset{\xi,s\geq 0}{\text{argmin}}\;\mathcal{L}^\text{ALM}_{\delta}(\xi,s;\lambda,\nu)+ \frac{\rho}{2}\left\|\xi-x\right\|_2^2,  \\\label{eq:lamupdate}
\lambda^+ & = \lambda+\frac{1}{\delta}\left(A x^+-b\right), \\\label{eq:nuupdate}
\nu^+ & = \nu+\frac{1}{\delta}\left(G x^+-h+s^+\right),
\end{align}
\end{subequations}
where superscript $^+$ denotes the iteration update.
Note that compared to the standard form for the original problem~\eqref{eq:qp_primal},~\eqref{eq:proximalM} does not add a penalty term for the slack variable as it does not contribute to the cost of the primal problem~\eqref{eq:qp_primal_slack}.

\section{Practical Algorithm}
\label{sec:pipm}

The method implemented in \texttt{PIQP}, following the framework proposed in~\cite{pougkakiotis2021}, deals with~\eqref{eq:PMM_sub_prob} by applying one iteration of the interior-point method, and uses the Mehrotra predictor-corrector
method~\cite{mehrotra1992} to update the primal-dual iterates. To detail the algorithm, we define the log-barrier function 
\begin{equation}
\Phi_\mu(s) = -\mu \sum_{i=1}^m \ln [s]_i,
\end{equation}
where $[s]_i$ denotes the $i$-th element of $s$, and $\mu>0$ is usually referred to as the barrier parameter. Replacing the constraint $s \geq 0$ in problem~\eqref{eq:PMM_sub_prob} with the penalty term $\sigma \cdot \Phi_\mu(s)$ in the objective function yields
\[
\min_{x,s}\;\;\mathcal{L}^\text{ALM}_{\delta_k}(x,s;\lambda_k,\nu_k)+\sigma_k\cdot \Phi_{\mu_k}(s)+ \frac{\rho_k}{2}\left\|x-\xi_k\right\|_2^2 
\]
at iteration $k$, where $(\xi_k,\lambda_k,\nu_k)$ defines the primal-dual iterates, and $\sigma_k\in(0,1]$ is the centering parameter in the predictor-corrector method discussed below. The first-order optimality conditions of the resulting unconstrained problem can then be represented as
\begin{subequations} 
\label{eq:PMM_IP_KKT}
\begin{align}
Px + c + \rho_k(x - \xi_k) + A^\top y + G^\top z & =0, \\
Ax - \delta_k (y-\lambda_k) - b & =0, \\
Gx - \delta_k (z-\nu_k) - h + s & =0, \\
s \circ z - \sigma_k\mu_k \mathbf{1}_m & =0, 
\end{align}
\end{subequations}
with 
auxiliary variables $y\in\mathbb R^p$ and $z\in\mathbb R^m$. Introducing auxiliary variables yields a sparser linear system of equations allowing highly efficient numerical routines. %

Applying Newton's method to solve~\eqref{eq:PMM_IP_KKT} results in the following linear equations 
\begin{equation} \label{eq:affine}
\underbrace{\begin{bmatrix}    
P + \rho_k I_n & A^\top & G^\top & 0 \\[0.12cm]
A & -\delta_k I_p & 0 & 0 \\[0.12cm]
G & 0 & -\delta_k I_m & I_m \\[0.12cm]
0 & 0 & S_k & Z_k
\end{bmatrix}}_{J(s_k,z_k)}
\underbrace{
\begin{bmatrix}
\Delta x_k \\[0.12cm]
\Delta y_k \\[0.12cm]
\Delta z_k \\[0.12cm]
\Delta s_k
\end{bmatrix}}_{\Delta \omega_k}
=
\underbrace{\begin{bmatrix}
r^x_k \\[0.12cm]
r^y_k \\[0.12cm]
r^z_k \\[0.12cm]
r^s_k
\end{bmatrix}}_{
r_k
}
\end{equation}
with initialization $(x_k,y_k,z_k,s_k)$ and 
\begin{equation*}
\begin{aligned}
r^x_k &= - (Px_k + c + \rho_k(x_k-\xi_k) + A^\top y_k + G^\top z_k), \\
r^y_k &= - (Ax_k + \delta_k(\lambda_k-y_k) - b), \\
r^z_k &= - (Gx_k + \delta_k(\nu_k - z_k) - h + s_k), \\
r^s_k &= -s_k \circ z_k + \sigma_k \mu_k\mathbf{1}_m
\end{aligned}
\end{equation*}
 at iteration $k$.
\begin{remark}
In the implementation, we can eliminate $\Delta s_k$. To this end, we compute the Nesterov-Todd scaling $W_k = Z_k^{-1} S_k$~\cite{nocedal2006} such that~\eqref{eq:affine} can be rewritten as
\[
\underbrace{\begin{bmatrix}
P + \rho_k I_n & A^\top & G^\top \\[0.12cm]
A & -\delta_k I_p & 0 \\[0.12cm]
G & 0 & -(W_k + \delta_k I_m)
\end{bmatrix}}_{\tilde J(s_k,z_k)}
\begin{bmatrix}
\Delta x_k \\[0.12cm]
\Delta y_k \\[0.12cm]
\Delta z_k
\end{bmatrix}
=\underbrace{\begin{bmatrix}
r^x_k \\[0.12cm]
r^y_k \\[0.12cm]
\bar{r}^z_k
\end{bmatrix}}_{\tilde r_k}
\]
with $\bar{r}^z_k = r^z_k - Z_k^{-1} r^s_k$. Here, $Z_k \in \mathbb{R}^{m \times m}$ is a diagonal matrix with $z_k$ on its diagonal.
Note that the slack direction $\Delta s_k$ 
can be reconstructed with $\Delta s_k = Z_k^{-1}(r^s_k - S_k \Delta z_k)$.
This also ensures that $\tilde J(s_k,z_k)$ is symmetric.
\end{remark}

Next, we present the three main steps of the Mehrotra predictor-corrector method:

\subsubsection{Prediction}
solve~\eqref{eq:affine} with $\mu_k=0$ and solution $\Delta \omega_k^a = (\Delta x_k^a, \Delta y_k^\text{a},\Delta z_k^a, \Delta s_k^\text{a})$
\subsubsection{Step Size and Centering Parameter} compute the primal and dual step sizes
\begin{align}\notag
\alpha^\text{a}_p =&\max\left\{\alpha\in[0,1]| s_k+\alpha \Delta s_k^\text{a}\geq (1-\tau)s_k\right\},\\\label{eq:stepsize}
\alpha^\text{a}_d =&\max\left\{\alpha\in[0,1]| z_k+\alpha \Delta z_k^\text{a}\geq (1-\tau)z_k\right\},
\end{align}
with the scaling parameter $\tau = 0.995$ chosen heuristically~\cite{pougkakiotis2021,nocedal2006} that ensures the iterates do not get too close to the boundary of the feasible set, and evaluate the centering parameter following~\cite{mehrotra1992}
\begin{equation}
\label{eq:parameter}
\sigma_k  =\max \{0, \min \{1, \eta_k\}\}^3
\end{equation}
with $\eta_k  =\left(\left(s_k+\alpha^a_p \Delta s^a_k\right)^\top\left(z_k+\alpha^a_d \Delta z^a_k\right)\right)/\left(\left(s_k^\top z_k\right) / m \right)$.
\subsubsection{Combined Correction and Centering} compute $\mu_k = \left(s_k^\top z_k\right) / m$ and solving~\eqref{eq:affine} with replacing $r_k^s$ by  
\begin{equation}
\label{eq:rsk}
r_k^s = \underbrace{- S_k z_k}_{\text{prediction}} \underbrace{- \Delta s^a_k \circ \Delta z^a_k}_{\text{correction}} + \underbrace{\sigma_k\cdot \mu_k\cdot \mathbf{1}_m}_{\text{centering}}
\end{equation}
yields solution $\Delta \omega_k^c = (\Delta x_k^\text{c}, \Delta y_k^\text{c},\Delta z_k^a, \Delta s_k^\text{c})$. Then, update $(x_{k+1},s_{k+1},y_{k+1},z_{k+1}):=$
\[
(x_{k},s_{k},y_{k},z_{k}) +  (\alpha^\text{c}_p\Delta x_{k}^\text{c},\alpha^\text{c}_p\Delta s_{k}^\text{c},\alpha^\text{c}_d\Delta y_{k}^\text{c},\alpha^\text{c}_d\Delta z_{k}^\text{c})
\]
with step sizes
\begin{align}\notag
\alpha^\text{c}_p =&\max\left\{\alpha\in[0,1]| s_k+\alpha \Delta s_k^\text{c}\geq (1-\tau)s_k\right\},\\\label{eq:centering_step_size}
\alpha^\text{c}_d =&\max\left\{\alpha\in[0,1]| z_k+\alpha \Delta z_k^\text{c}\geq (1-\tau)z_k\right\}.
\end{align}

Based on the aforementioned discussion, the primal-dual iterate $(\xi_k,\lambda_k,\nu_k)$ is optionally updated as outlined in Algorithm~\ref{alg:PEU} proposed by~\cite[Section 5.1.4]{pougkakiotis2021}. Here, we introduce the primal-dual residual with respect to the $k$-th iteration
\[
\begin{aligned}
p_k =& \left\|\begin{bmatrix} Ax_k - b \\ Gx_k - h + s_k \end{bmatrix}\right\|_\infty,\\
d_k =& \left\|Px_k + c +A^\top y_k + G^\top z_k\right\|_\infty.
\end{aligned}
\]
Moreover, $\delta$ and $\rho$ are limited by $\underline{\delta}$ and $\underline{\rho}$ for numerical stability.

Now, we can summarize a predictor-corrector-based practical computational framework for the interior-point proximal 
method of multipliers to solve~\eqref{eq:qp_primal_slack} in Algorithm~\ref{alg:ProxIP}.
\begin{remark}
Algorithm~\ref{alg:ProxIP} is a practical variant of the standard interior-point proximal method of multipliers as presented in~\cite[Algorithm~1]{pougkakiotis2021}, which substitutes the correction step and regularization update in~\cite[Algorithm~1]{pougkakiotis2021} by the Mehrotra's predictor-corrector method and Algorithm~\ref{alg:PEU}, respectively. Note that Algorithm~\ref{alg:PEU} is a numerical heuristic to update the primal-dual iterates $(\xi_k,\lambda_k,\nu_k)$ such that the convergence analysis proposed in~\cite[Section~3]{pougkakiotis2021} cannot be rigorously established step by step. However, this heuristic leads to a more reliable and effective numerical convergence, although compared to the update in~\cite[Algorithm~1]{pougkakiotis2021} without theoretical guarantees.  
Analyzing the theoretical convergence guarantee of practical Algorithm~\ref{alg:ProxIP} is beyond the scope of this paper and will be investigated in our future work.    
\end{remark}
\vspace{-0.3cm}
\begin{algorithm}[htbp!]
\caption{Interior-Point Proximal Method of Multipliers for Convex Quadratic Programming}
\label{alg:ProxIP}
\small
\textbf{Initialization:} choose $(\xi_0,s_0,\lambda_0,\nu_0)$, set $\delta_0,\rho_0>0$.
\begin{algorithmic}[1]
\For{$k=0,1,...$}
\State \hspace{-3mm}Set $(x_k,  y_k,z_k)=(\xi_k,\lambda_k,\nu_k)$;
\State \hspace{-3mm}Solve~\eqref{eq:affine} with $r_k^s = -S_kz_k$ for 
$\Delta \omega_k^\text{a}$;
\State \hspace{-3mm}Compute $(\alpha_p^\text{a},\alpha_d^\text{a})$ by~\eqref{eq:stepsize};
\State \hspace{-3mm}Update $\mu_k=\left(s_k^\top z_k\right)/m$ and $\sigma_k$ by~\eqref{eq:parameter};
\State \hspace{-3mm}Solve~\eqref{eq:affine} with $r_k^s$ by~\eqref{eq:rsk} for $\Delta\omega_k^\text{c}$;
\State \hspace{-3mm}Compute $(\alpha^\text{c}_p,\alpha^\text{c}_d)$ by~\eqref{eq:centering_step_size} and update 
\[
\begin{aligned}
(&x_{k+1},s_{k+1},y_{k+1},z_{k+1}) \leftarrow    \\
&\quad (x_{k},s_{k},y_{k},z_{k}) +  (\alpha^\text{c}_p\Delta \xi_{k}^\text{c},\alpha^\text{c}_p\Delta s_{k}^\text{c},\alpha^\text{c}_d\Delta y_{k}^\text{c},\alpha^\text{c}_d\Delta z_{k}^\text{c});
\end{aligned}
\]
\State \hspace{-3mm}Run Alg.~\ref{alg:PEU} to get $(\xi_{k+1},\lambda_{k+1},\nu_{k+1})$, $(\delta_{k+1},\rho_{k+1})$. 
\EndFor
\end{algorithmic}
\end{algorithm}
\vspace{-0.5cm}
\begin{algorithm}[htbp!]
\caption{Penalty and Estimate Updates}
\label{alg:PEU}
\small
\textbf{Input:} $r=\left|s_{k}^\top z_{k} - s_{k+1}^\top z_{k+1}\right|/s_{k}^\top z_{k}$.\\
\textbf{Output:} $(\xi_{k+1},\lambda_{k+1},\nu_{k+1})$ and $(\delta_{k+1},\rho_{k+1})$. 
\begin{algorithmic}[1]
\If{$p_{k+1}\leq 0.95 \cdot p_k$}
\State 
$(\lambda_{k+1},\nu_{k+1})\leftarrow (y_{k+1},z_{k+1})$, $\delta_{k+1} \leftarrow (1-r) \delta_k$
\Else
\State $(\lambda_{k+1},\nu_{k+1})\leftarrow (\lambda_{k},\nu_{k})$, $\delta_{k+1} \leftarrow (1-r/3)\delta_k$
\EndIf
\If{$d_{k+1}\leq 0.95 \cdot d_k$}
\State  $\xi_{k+1}\leftarrow x_{k+1}$, $\rho_{k+1} \leftarrow (1-r) \rho_k$
\Else
\State $\xi_{k+1}\leftarrow \xi_k$, $\rho_{k+1} \leftarrow (1-r/3)\rho_k$
\EndIf
\State $\delta_{k+1} \leftarrow \max\{\delta_{k+1},\underline{\delta}\}$, $\rho_{k+1} \leftarrow \max\{\rho_{k+1},\underline{\rho}\}$%
\end{algorithmic}
\end{algorithm}
\vspace{-0.5cm}
\section{Numerical Implementation} \label{sec:imp}

\subsection{Initialization}
We initialize the primal and dual variables using the standard method proposed in \cite{andersen2011} by minimizing the unconstrained optimization problem
\begin{equation*}
\begin{aligned}
\min_{\xi_0,\tilde{s}_0} \quad & \mathcal{L}^\text{ALM}_{\delta_0}(\xi_0,\tilde{s}_0;0,0)+ \frac{\rho_0}{2}\left\|\xi_0\right\|_2^2 + \frac{1}{2} \|\tilde{s}_0\|_2^2,
\end{aligned}
\end{equation*}
which can be posed as the solution to the linear system of equations
\begin{equation}
\begin{bmatrix}
P + \rho_0 I_n & A^\top & G^\top \\[0.12cm]
A & -\delta_0 I_p & 0 \\[0.12cm]
G & 0 & -(1 + \delta_0) I_m
\end{bmatrix}
\begin{bmatrix}
\xi_0 \\[0.12cm]
\lambda_0 \\[0.12cm]
\tilde{\nu}_0
\end{bmatrix}
=
\begin{bmatrix}
-c \\[0.12cm]
b \\[0.12cm]
h
\end{bmatrix}
\end{equation}
with potentially negative slack variable $\tilde{s}_0 = -\tilde{\nu}_0$. Note that the structure is the same as for $\tilde J(\tilde{s}_0,\tilde{\nu}_0)$; hence, we can reuse the symbolic factorization, resulting in lower computational cost.

To guarantee that $s_0$ and $\nu_0$ are in the non-negative orthant with sufficient magnitude, we calculate the conservative step sizes as
\begin{equation*}
\begin{aligned}
\Delta \tilde{s}_0 &= \max \{0, -1.5 \cdot \min(\tilde{s}_0)\}, \\
\Delta \tilde{\nu}_0 &= \max \{0, -1.5 \cdot \min(\tilde{\nu}_0)\},
\end{aligned}
\end{equation*}
and similar to \cite{pougkakiotis2021}, we shift initial solutions further away from the barrier based on the normalized complementarity violation
\begin{equation*}
\begin{aligned}
\Delta s_0 &= \Delta \tilde{s}_0 + 0.5 \cdot \frac{(\tilde{s}_0 + \Delta \tilde{s}_0)^\top(\tilde{\nu}_0 + \Delta \tilde{\nu}_0)}{\sum_{i=0}^m ([\tilde{\nu}_0]_i + [\Delta \tilde{\nu}_0]_i)}, \\
\Delta \nu_0 &= \Delta \tilde{\nu}_0 + 0.5 \cdot \frac{(\tilde{s}_0 + \Delta \tilde{s}_0)^\top(\tilde{\nu}_0 + \Delta \tilde{\nu}_0)}{\sum_{i=0}^m ([\tilde{s}_0]_i + [\Delta \tilde{s}_0]_i)},
\end{aligned}
\end{equation*}
resulting in the initial values for slack and inequality Lagrange multipliers 
$s_0 = \tilde{s}_0 + \Delta s_0$, $\nu_0 = \tilde{\nu}_0 + \Delta \nu_0$

\subsection{Termination Criteria}

We adopt the same terminal criteria for convergence as in \texttt{SCS v3.0} \cite{donoghue2021}. More specifically, \texttt{PIQP} terminates when it finds primal variables $x\in \mathbb{R}^n$, $s \in \mathbb{R}^m$, and dual variables $y \in \mathbb{R}^p$, $z \in \mathbb{R}^m$ which satisfy the conditions
\begin{subequations}
\begin{align}
&\left\|\begin{bmatrix} Ax - b \\ Gx - h + s \end{bmatrix}\right\|_{\infty} \leq \epsilon_{\mathrm{abs}} \label{eq:term_primal} \\
& \;\;\qquad +\epsilon_{\mathrm{rel}} \max \left(\|A x\|_{\infty},\|b\|_{\infty},\|Gx\|_{\infty},\|h\|_{\infty},\|s\|_{\infty}\right), \nonumber \\[0.12cm]
&\left\|P x+A^\top y + G^\top z +c\right\|_{\infty} \leq \epsilon_{\mathrm{abs}} \label{eq:term_dual} \\
& \;\;\qquad +\epsilon_{\mathrm{rel}} \max \left(\|P x\|_{\infty}, \left\|A^\top y\right\|_{\infty},\left\|G^\top z\right\|_{\infty},\|c\|_{\infty}\right), \nonumber \\[0.12cm]
&\left|x^\top P x+c^\top x+b^\top y + h^\top z\right| \leq \epsilon_{\mathrm{abs}} \label{eq:term_gap} \\
& \;\;\qquad +\epsilon_{\mathrm{rel}} \max \left(\left|x^\top P x\right|,\left|c^\top x\right|,\left|b^\top y\right|,\left|h^\top z\right|\right), \nonumber
\end{align}
\end{subequations}
where $\epsilon_{\mathrm{abs}} > 0$ and $\epsilon_{\mathrm{rel}} \geq 0$ are the user defined absolute and relative accuracies. Condition \eqref{eq:term_primal} corresponds to the primal feasibility, and \eqref{eq:term_dual} to the dual feasibility, which is common in most solvers like \texttt{OSQP}~\cite{stellato2020} or \texttt{qpSWIFT}~\cite{pandala2019}. The condition on the duality gab \eqref{eq:term_gap} is less commonly checked, but if neglected, it can result in poor solution quality. The Maros-Mészáros problem set, for example, includes problems that are solved inaccurately without the criteria on the duality gap as discussed in \cite[Section~7.2]{donoghue2021}.

\subsection{Sparse Pivot-Free LDL Factorization} \label{subsec:factorization}

The solution of the KKT system \eqref{eq:affine} constitutes the most computationally expensive step in any interior-point method. In our work, we have chosen to employ a direct method that is particularly well-suited for use in embedded applications. Specifically, we have opted for a modified version of Tim Davis' sparse pivot-free LDL factorization \cite{davis2005} with approximate minimum degree (AMD) ordering \cite{amestoy2004}, resulting in the factorization
\begin{equation*}
\Gamma K\Gamma ^\top = LDL^\top,
\end{equation*}
where $\Gamma$ is the permutation matrix of the AMD ordering reducing fill-in of the lower triangular matrix $L$, and $D$ is a diagonal matrix.

To ensure that the LDL factorization of any symmetric permutation of the KKT matrix exists, it is sufficient if $K$ is quasi-definite \cite{vanderbei1995}. Although adding terms to the diagonal of \eqref{eq:affine} through the proximal method of multipliers typically ensures that $K$ is almost certainly quasi-definite, there may be rare instances when it is not. In such cases, we slightly perturb the regularization parameters and retry the factorization. Thanks to this approach, there is no need to resort to techniques such as dynamic regularization with subsequent iterative refinement, which are employed in \texttt{ECOS}~\cite{domahidi2013}, resulting in less computational overhead.

The LDL factorization process consists of two phases: symbolic and numeric. During the symbolic phase, the elimination tree and the fill-in pattern of the lower triangular matrix $L$ are determined, which provides the necessary information regarding the required memory. In the numeric phase, the factorization of $K$ is performed using the previously calculated elimination tree, which fills in both $L$ and $D$. Since the structure of $\tilde J(s_k,z_k)$ remains constant, we can reuse the elimination tree and fill-in pattern for every subsequent solve. As a result, symbolic computation can be avoided, and all previously allocated memory can be reused. This provides a significant advantage in terms of computational efficiency, particularly in the context of optimal control, where the structure of the sequence of solved problems stays constant.

\subsection{Preconditioning}

Preconditioning is a well-established technique for reducing the number of iterations for first-order methods by decreasing the condition number of the KKT system \cite[Chapter~5]{nocedal2006}. Although commonly used for first-order methods, preconditioners can also be beneficial for interior-point methods by improving numerical stability and convergence. In our implementation, we adopt the Ruiz equilibration as our default preconditioner \cite{ruiz2001}. This technique scales the problem data diagonally to reduce the conditioning number of the unregularized KKT conditions \eqref{eq:affine} while being relatively cheap to compute.

\begin{remark}
With the implementation of the Ruiz equilibration technique, we observed a notable improvement, achieving an average speedup of $22\%$ on the Maros-Mészáros problem set. Before preconditioning, one less problem could be solved.
\end{remark}

\subsection{Interface and Memory Allocation} \label{subsec:interface}

In addition to the QP formulation \eqref{eq:qp_primal} discussed in the paper, our implementation includes an interface for incorporating box constraints of the form $l \leq x \leq u$ with $l, u \in \mathbb{R}^n$. By exploiting this underlying structure, we are able to achieve computational speedups when factoring the corresponding KKT matrix.

Our QP solver, \texttt{PIQP}, leverages the high performance of \texttt{C/C++} and the efficient vectorized matrix/vector operations provided by the popular \texttt{Eigen3} library \cite{eigenweb}. In addition to the \texttt{C/C++} interface, we provide a \texttt{Python}, \texttt{Matlab}, and an \texttt{R} interface, making it easy to integrate \texttt{PIQP} into a wide range of computational pipelines. The code structure and Python interface have been adopted from ProxSuite~\cite{bambade2022}.

The interface of \texttt{PIQP} is comprised of three main routines: \texttt{setup}, \texttt{update}, and \texttt{solve}. During the \texttt{setup} routine, all necessary memory structures are set up, the problem data is preconditioned, and the symbolic factorization, including the AMD ordering, is computed. The \texttt{update} routine is designed for updating problem data in subsequent solves, reusing the existing memory and factorization, assuming the sparsity structure remains unchanged. Similar to the \texttt{setup}, the updated problem data gets preconditioned. Upon executing the algorithm, the \texttt{solve} routine returns the solution, as well as the status of the optimization.

In predictive optimal control applications, it is common to solve the same problem repeatedly without any changes in structure. By utilizing the \texttt{update} routine, the problem data can be updated while reusing the symbolic factorization and already allocated memory. This significantly speeds up subsequent solutions.

Similar to solvers like \texttt{OSQP}~\cite{stellato2020}, \texttt{ECOS}~\cite{domahidi2013}, and \texttt{PROXQP}~\cite{bambade2022}, dynamic memory is only allocated during the \texttt{setup} routine, ensuring that the \texttt{update} and \texttt{solve} functions remain \texttt{malloc}-free. This feature is particularly important for embedded systems, where memory should be allocated statically or at least only once, without any subsequent dynamic allocations.

\section{Numerical Examples}
\label{sec:num}

To demonstrate the performance of \texttt{PIQP}, we benchmark it against the open-source toolkits \texttt{OSQP}~\cite{stellato2020}, \texttt{SCS}~\cite{donoghue2016}, and \texttt{PROXQP}~\cite{bambade2022}, as well as commercial solvers \texttt{Gurobi}~\cite{gurobi} and \texttt{Mosek}~\cite{mosek}. Note that we did not compare our solver with the dense solver \texttt{qpOASES}~\cite{ferreau2014qpoases} or conic programming solvers such as \texttt{ECOS}~\cite{domahidi2013}, as they do not support quadratic objectives natively.

To make the comparison fair, we enforced the same termination criteria across all solvers. Specifically, all solvers verify the primal feasibility \eqref{eq:term_primal} and dual feasibility \eqref{eq:term_dual}. However, the duality gap \eqref{eq:term_gap} is only checked by \texttt{PIQP}, \texttt{SCS}, and \texttt{PROXQP}. As a result, other solvers may report an optimal solution that fails to satisfy the duality gap condition and may require more time or even fail to find a solution that meets our termination criteria. Thus, it is important to note that the reported results for these solvers may be overly optimistic. For example, \texttt{Mosek} would fail on $75 \%$ of problems if we check the termination criteria of the solution.

All benchmarks are run on a workstation with an AMD Ryzen
Threadripper 3990X 4.3 GHz CPU. \texttt{Gurobi} and \texttt{Mosek} were limited to a single thread, and all solvers were subject to a $1000$ second time limit. We use an adapted benchmark framework from \texttt{OSQP}~\cite{stellato2020}.

\subsection{Maros-Mészáros problems} 
\label{subsec:maros}

We consider the standard Maros-Mészáros problem set, comprised out of $138$ \textit{hard} QP problems \cite{maros1999}. Most problems are very sparse, and a certain subset is numerically extremely ill-conditioned.

Our implementation is based on the \texttt{Python} interface of all solvers. Moreover, we only use the internally measured times reported by the solver. This includes setup and solution time.

We conduct two sets of experiments to evaluate the performance of our QP solver. In the first scenario, we aim to find solutions with low accuracy, setting $\epsilon_{\mathrm{abs}} = 10^{-3}$ and $\epsilon_{\mathrm{rel}} = 10^{-4}$. In the second scenario, we target highly accurate solutions, setting $\epsilon_{\mathrm{abs}} = 10^{-8}$ and $\epsilon_{\mathrm{rel}} = 10^{-9}$.

Table~\ref{table:maros_meszaros_failure_rates} summarizes the failure rates of the solvers for the different accuracy settings. Compared to the other solvers, our solver, \texttt{PIQP}, demonstrates remarkable robustness, successfully solving all problems with low accuracy settings and failing to solve only a single problem with high accuracy settings.
\vspace{-0.5cm}
\begin{table}[htbp!]
\centering
\caption{Failure rates on the Maros-Mészáros problem set}
\renewcommand{\arraystretch}{1.5}
\tiny
\begin{tabular}{lcccccc} 
\hline
& \texttt{PIQP} &
\texttt{OSQP} & 
\texttt{SCS} & 
\texttt{PROXQP} & \texttt{GUROBI} & 
\texttt{MOSEK} \\\hline
Low Accuracy &  $0.00 \%$ & $4.35 \%$ & $4.35 \%$ & $12.32 \%$ & $6.52 \%$  & $5.07 \%$\\
\hline
High Accuracy &  $0.72 \%$  & $31.16 \%$& $23.19 \%$ & $23.19 \%$ & $6.52 \%$  & $10.14 \%$ \\
\hline
\end{tabular}
\label{table:maros_meszaros_failure_rates}
\end{table}

\begin{figure}[htbp!]
\centering
\includegraphics[width=0.48\textwidth]{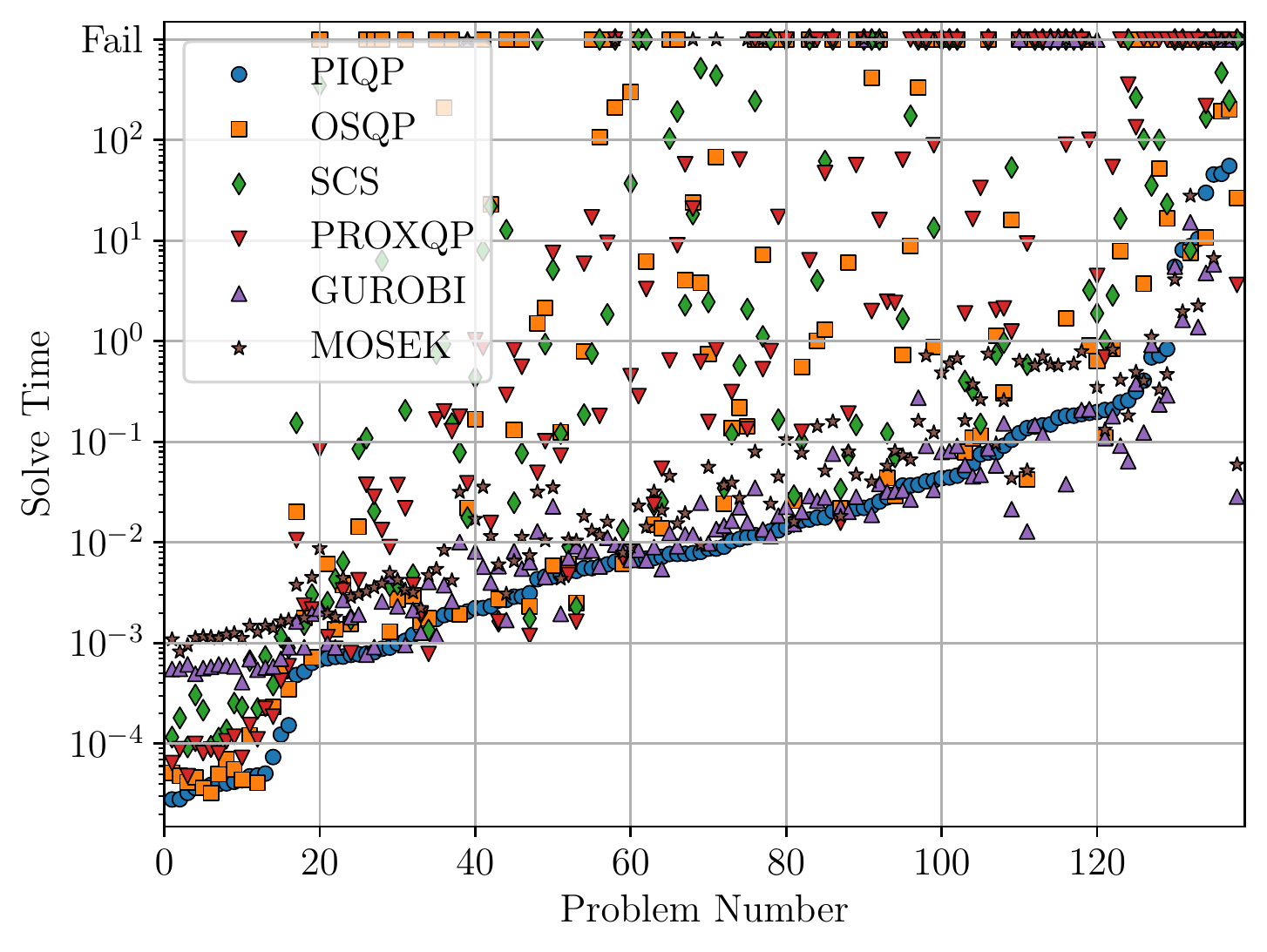}
\vspace{-0.5cm}
\caption{Solve times for high accuracy settings for individual problems in the Maros-Mészáros problem set, ordered by \texttt{PIQP} solve time. The lowest point in each column corresponds to the fastest solver.}
\label{fig:maros_meszaros_problems_high_accuracy_time}
\end{figure}

We provide the individual solve times for each problem in the Maros-Meszaros problem set for high accuracy settings in Figure~\ref{fig:maros_meszaros_problems_high_accuracy_time}. The problems are sorted by \texttt{PIQP} solve time, highlighting the speed of our solver relative to the others. Notably, our solver is almost always the fastest for most problems. These results highlight the advantages of interior-point methods in high-accuracy settings, where they outperform most first-order methods.

\begin{figure}[htbp!]
\centering
\includegraphics[width=0.48\textwidth]{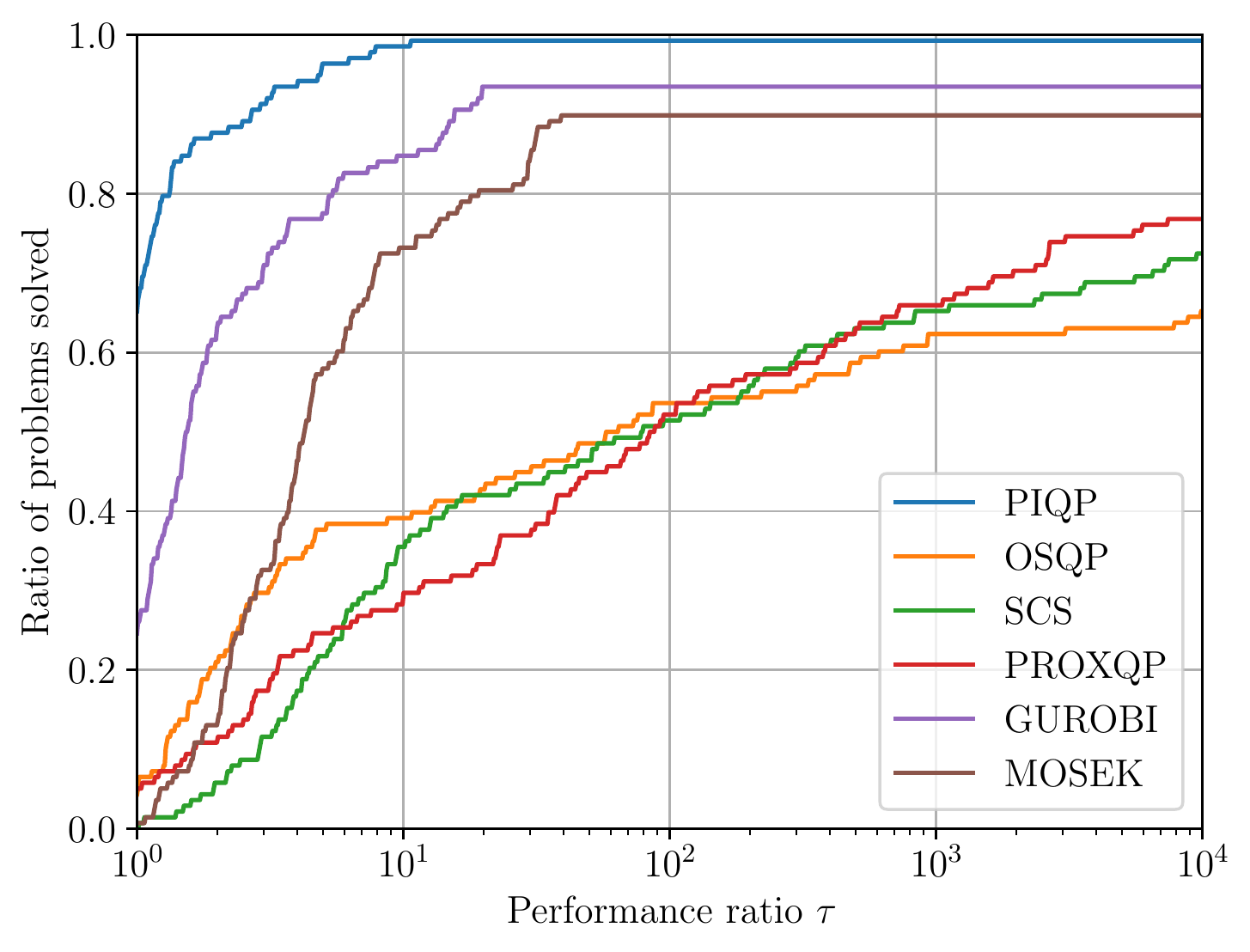}
\vspace{-0.5cm}
\caption{Performance profiles on the Maros-Mészáros problem set for high accuracy settings.}
\label{fig:maros_meszaros_problems_high}
\end{figure}

We include the Dolan-Moré performance profiles \cite{dolan2002} for the high accuracy settings in Figure~\ref{fig:maros_meszaros_problems_high}. The x-axis of the graphs shows the solve time normalized by the fastest solver, with a performance ratio of one indicating that the solver was the fastest and a performance ratio of ten indicating that the solver was ten times slower than the fastest solver for a particular problem. The y-axis of the graphs shows the ratio of problems solved. For instance, \texttt{PIQP} was the fastest solver for 65\% of the problems and took at most five times longer for 96\% of the problems.

\bibliographystyle{IEEEtran}
\bibliography{refs}

\begin{thebibliography}{10}
\providecommand{\url}[1]{#1}
\csname url@samestyle\endcsname
\providecommand{\newblock}{\relax}
\providecommand{\bibinfo}[2]{#2}
\providecommand{\BIBentrySTDinterwordspacing}{\spaceskip=0pt\relax}
\providecommand{\BIBentryALTinterwordstretchfactor}{4}
\providecommand{\BIBentryALTinterwordspacing}{\spaceskip=\fontdimen2\font plus
\BIBentryALTinterwordstretchfactor\fontdimen3\font minus
  \fontdimen4\font\relax}
\providecommand{\BIBforeignlanguage}[2]{{%
\expandafter\ifx\csname l@#1\endcsname\relax
\typeout{** WARNING: IEEEtran.bst: No hyphenation pattern has been}%
\typeout{** loaded for the language `#1'. Using the pattern for}%
\typeout{** the default language instead.}%
\else
\language=\csname l@#1\endcsname
\fi
#2}}
\providecommand{\BIBdecl}{\relax}
\BIBdecl

\bibitem{rubinstein2002markowitz}
M.~Rubinstein, ``Markowitz's" portfolio selection": A fifty-year
  retrospective,'' \emph{The Journal of {F}inance}, vol.~57, no.~3, pp.
  1041--1045, 2002.

\bibitem{rawlings2017model}
J.~B. Rawlings, D.~Q. Mayne, and M.~Diehl, \emph{Model {P}redictive {C}ontrol:
  {T}heory, {C}omputation, and {D}esign}.\hskip 1em plus 0.5em minus
  0.4em\relax Nob Hill Publishing Madison, WI, 2017.

\bibitem{allgower1999nonlinear}
F.~Allg{\"o}wer, T.~A. Badgwell, J.~S. Qin, J.~B. Rawlings, and S.~J. Wright,
  ``Nonlinear predictive control and moving horizon estimation—an
  introductory overview,'' \emph{Advances in control: Highlights of ECC’99},
  pp. 391--449, 1999.

\bibitem{zhu2018blended}
Y.~Zhu, R.~Bridson, and D.~M. Kaufman, ``Blended cured quasi-{N}ewton for
  distortion optimization,'' \emph{ACM Transactions on Graphics}, vol.~37,
  no.~4, pp. 1--14, 2018.

\bibitem{boggs1995sequential}
P.~T. Boggs and J.~W. Tolle, ``Sequential quadratic programming,'' \emph{Acta
  {N}umerica}, vol.~4, pp. 1--51, 1995.

\bibitem{fletcher1998numerical}
R.~Fletcher and S.~Leyffer, ``Numerical experience with lower bounds for {MIQP}
  branch-and-bound,'' \emph{SIAM Journal on Optimization}, vol.~8, no.~2, pp.
  604--616, 1998.

\bibitem{bemporad2015quadratic}
A.~Bemporad, ``A quadratic programming algorithm based on nonnegative least
  squares with applications to embedded model predictive control,'' \emph{IEEE
  Transactions on Automatic Control}, vol.~61, no.~4, pp. 1111--1116, 2015.

\bibitem{cimini2019complexity}
G.~Cimini and A.~Bemporad, ``Complexity and convergence certification of a
  block principal pivoting method for box-constrained quadratic programs,''
  \emph{Automatica}, vol. 100, pp. 29--37, 2019.

\bibitem{ferreau2014qpoases}
H.~J. Ferreau, C.~Kirches, A.~Potschka, H.~G. Bock, and M.~Diehl, ``{qpOASES:
  A} parametric active-set algorithm for quadratic programming,''
  \emph{Mathematical Programming Computation}, vol.~6, pp. 327--363, 2014.

\bibitem{wright1997}
S.~J. Wright, \emph{Primal-{D}ual {I}nterior-{P}oint {M}ethods}.\hskip 1em plus
  0.5em minus 0.4em\relax SIAM, 1997.

\bibitem{domahidi2013}
A.~Domahidi, E.~Chu, and S.~Boyd, ``{ECOS: An SOCP} solver for embedded
  systems,'' in \emph{European Control Conference}, 2013, pp. 3071--3076.

\bibitem{frison2020hpipm}
G.~Frison and M.~Diehl, ``{HPIPM}: {A} high-performance quadratic programming
  framework for model predictive control,'' \emph{IFAC-PapersOnLine}, vol.~53,
  no.~2, pp. 6563--6569, 2020.

\bibitem{goldstein2014fast}
T.~Goldstein, B.~O'Donoghue, S.~Setzer, and R.~Baraniuk, ``Fast alternating
  direction optimization methods,'' \emph{SIAM Journal on Imaging Sciences},
  vol.~7, no.~3, pp. 1588--1623, 2014.

\bibitem{patrinos2013accelerated}
P.~Patrinos and A.~Bemporad, ``An accelerated dual gradient-projection
  algorithm for embedded linear model predictive control,'' \emph{IEEE
  Transactions on Automatic Control}, vol.~59, no.~1, pp. 18--33, 2013.

\bibitem{patrinos2015dual}
P.~Patrinos, A.~Guiggiani, and A.~Bemporad, ``A dual gradient-projection
  algorithm for model predictive control in fixed-point arithmetic,''
  \emph{Automatica}, vol.~55, pp. 226--235, 2015.

\bibitem{stellato2020}
B.~Stellato, G.~Banjac, P.~Goulart, A.~Bemporad, and S.~Boyd, ``{OSQP: A}n
  operator splitting solver for quadratic programs,'' \emph{Mathematical
  Programming Computation}, vol.~12, no.~4, pp. 637--672, 2020.

\bibitem{frasch2015parallel}
J.~V. Frasch, S.~Sager, and M.~Diehl, ``A parallel quadratic programming method
  for dynamic optimization problems,'' \emph{Mathematical Programming
  Computation}, vol.~7, pp. 289--329, 2015.

\bibitem{patrinos2011global}
P.~Patrinos, P.~Sopasakis, and H.~Sarimveis, ``A global piecewise smooth
  {N}ewton method for fast large-scale model predictive control,''
  \emph{Automatica}, vol.~47, no.~9, pp. 2016--2022, 2011.

\bibitem{bambade2022}
A.~Bambade, S.~El-Kazdadi, A.~Taylor, and J.~Carpentier, ``{PROX-QP}: yet
  another quadratic programming solver for robotics and beyond,'' in
  \emph{{Robotics: Science and Systems}}, 2022.

\bibitem{liao2018regularized}
D.~Liao-McPherson, M.~Huang, and I.~Kolmanovsky, ``A regularized and smoothed
  {F}ischer--{B}urmeister method for quadratic programming with applications to
  model predictive control,'' \emph{IEEE Transactions on Automatic Control},
  vol.~64, no.~7, pp. 2937--2944, 2018.

\bibitem{bemporad2017numerically}
A.~Bemporad, ``A numerically stable solver for positive semidefinite quadratic
  programs based on nonnegative least squares,'' \emph{IEEE Transactions on
  Automatic Control}, vol.~63, no.~2, pp. 525--531, 2017.

\bibitem{hermans2022qpalm}
B.~Hermans, A.~Themelis, and P.~Patrinos, ``{QPALM}: {A} proximal augmented
  {L}agrangian method for nonconvex quadratic programs,'' \emph{Mathematical
  Programming Computation}, vol.~14, no.~3, pp. 497--541, 2022.

\bibitem{liao2020fbstab}
D.~Liao-McPherson and I.~Kolmanovsky, ``{FBstab}: A proximally stabilized
  semismooth algorithm for convex quadratic programming,'' \emph{Automatica},
  vol. 113, p. 108801, 2020.

\bibitem{pougkakiotis2021}
S.~Pougkakiotis and J.~Gondzio, ``An interior point-proximal method of
  multipliers for convex quadratic programming,'' \emph{Computational
  Optimization and Applications}, vol.~78, no.~2, pp. 307--351, 2021.

\bibitem{maros1999}
I.~Maros and C.~Mészáros, ``A repository of convex quadratic programming
  problems,'' \emph{Optimization Methods and Software}, vol.~11, no. 1-4, pp.
  671--681, 1999.

\bibitem{rockafellar1976}
R.~T. Rockafellar, ``Augmented {L}agrangians and applications of the proximal
  point algorithm in convex programming,'' \emph{Mathematics of Operations
  Research}, vol.~1, no.~2, pp. 97--116, 1976.

\bibitem{mehrotra1992}
S.~Mehrotra, ``On the implementation of a primal-dual interior point method,''
  \emph{SIAM Journal on Optimization}, vol.~2, no.~4, pp. 575--601, 1992.

\bibitem{nocedal2006}
J.~Nocedal and S.~J. Wright, \emph{Numerical Optimization}.\hskip 1em plus
  0.5em minus 0.4em\relax Springer, 2006.

\bibitem{andersen2011}
M.~Andersen, J.~Dahl, Z.~Liu, and L.~Vandenberghe, ``Interior-point methods for
  large-scale cone programming,'' in \emph{{Optimization for Machine
  Learning}}.\hskip 1em plus 0.5em minus 0.4em\relax The MIT Press, 2011.

\bibitem{donoghue2021}
B.~O'Donoghue, ``Operator splitting for a homogeneous embedding of the linear
  complementarity problem,'' \emph{SIAM Journal on Optimization}, vol.~31,
  no.~3, pp. 1999--2023, 2021.

\bibitem{pandala2019}
A.~G. Pandala, Y.~Ding, and H.-W. Park, ``{qpSWIFT: A} real-time sparse
  quadratic program solver for robotic applications,'' \emph{IEEE Robotics and
  Automation Letters}, vol.~4, no.~4, pp. 3355--3362, 2019.

\bibitem{davis2005}
T.~A. Davis, ``{Algorithm 849: A concise sparse Cholesky factorization
  package},'' \emph{ACM Transactions on Mathematical Software}, vol.~31, no.~4,
  pp. 587--591, 2005.

\bibitem{amestoy2004}
P.~R. Amestoy, T.~A. Davis, and I.~S. Duff, ``{Algorithm 837: AMD, an
  approximate minimum degree ordering algorithm},'' \emph{ACM Transactions on
  Mathematical Software}, vol.~30, no.~3, pp. 381--388, 2004.

\bibitem{vanderbei1995}
R.~J. Vanderbei, ``Symmetric quasidefinite matrices,'' \emph{SIAM Journal on
  Optimization}, vol.~5, no.~1, pp. 100--113, 1995.

\bibitem{ruiz2001}
R.~Daniel, ``A scaling algorithm to equilibrate both rows and columns norms in
  matrices,'' Rutherford Appleton Laboratory, Tech. Rep., 2001.

\bibitem{eigenweb}
G.~Guennebaud, B.~Jacob \emph{et~al.}, ``Eigen v3,''
  http://eigen.tuxfamily.org, 2010.

\bibitem{donoghue2016}
B.~O'Donoghue, E.~Chu, N.~Parikh, and S.~Boyd, ``Conic optimization via
  operator splitting and homogeneous self-dual embedding,'' \emph{Journal of
  Optimization Theory and Applications}, vol. 169, no.~3, pp. 1042--1068, 2016.

\bibitem{gurobi}
{Gurobi Optimization, LLC}, ``{Gurobi Optimizer Reference Manual},'' 2023.

\bibitem{mosek}
M.~ApS, \emph{MOSEK Optimizer API for Python 10.0.39}, 2023.

\bibitem{dolan2002}
E.~D. Dolan and J.~J. Mor{\'e}, ``Benchmarking optimization software with
  performance profiles,'' \emph{Mathematical Programming}, vol.~91, no.~2, pp.
  201--213, 2002.

\end{thebibliography}

\end{document}